\documentclass[12pt]{article}
\input{amssym.def}
\input{amssym.tex}
\newcommand{\ov}{\overline}

\newcommand{\ovp}{\ov p}

\newcommand{\dom}{\mbox{\rm dom}}
\newcommand{\qed}{{$\dashv$}}

\newcommand{\pregap}{\{ (a_i\mid i \in \omega_1), (b_j \mid j \in \omega_1)\} }

\newcommand{\IJpregap}{\{ (a_i\mid i \in I), (b_j \mid j \in J)\} }

\newcommand{\forces}{{\: \Vdash\:}}
\newcommand{\forcesQ}{{\: \Vdash_Q\:}}
\newcommand{\forcesR}{{\: \Vdash_R\:}}
\newcommand{\forcesP}{{\: \Vdash_P\:}}

\newcommand{\height}{\mbox{\rm height}}
\newcommand{\cf}{\mbox{\rm cf}}

\def\newtheorems{\newtheorem{theorem}{Theorem}[section]
                 \newtheorem{cor}[theorem]{Corollary}
                 \newtheorem{prop}[theorem]{Proposition}
                 \newtheorem{lemma}[theorem]{Lemma}

                 \newtheorem{definition}[theorem]{Definition}

                 }

\newtheorems

\title{Ladder gaps over stationary sets}
\author
{Uri Abraham\\
Departments of Mathematics and Computer Science, \\
Ben-Gurion University, Be\'{e}r-Sheva, Israel\\
Saharon Shelah 
\thanks{The author would like to thank the Israel
Science Foundation, founded by the Israel Academy of Science and
Humanities. Publication \# 598. }
\\
Institute of Mathematics\\
The Hebrew University, Jerusalem, Israel
}
\begin{document}
\bibliographystyle{plain}
\maketitle
\begin{abstract}
For a stationary set $S\subseteq \omega_1$ and a ladder system $C$ over
$S$,
a new type of gaps called $C$-Hausdorff
is introduced and investigated.
 We describe a forcing model of ZFC in which, for some stationary 
 set $S$, for every ladder $C$ over $S$, every gap
contains a subgap that is $C$-Hausdorff. But for every ladder $E$ over
$\omega_1\setminus S$ there exists a gap with
no subgap that is $E$-Hausdorff.

A new type of chain condition, called polarized chain
condition, is introduced. We prove that the iteration with finite
support of polarized c.c.c posets is again a polarized c.c.c poset.
\end{abstract}

\section{Introduction}
We first review some notations and definitions related to 
Hausdorff gaps. In fact 
we follow here the terminology given by M. Scheepers in his monograph
\cite{Scheepers} on Hausdorff gaps, but
 since we restrict ourselves to $(\omega_1,\omega^*_1)$ gaps 
our nomenclature is
somewhat simpler. The collection of all infinite subsets of $\omega$ is
denoted $[\omega ]^\omega$, and for $a,b \in [\omega ]^\omega$, $a
\subseteq^* b$ means that $a\setminus b$ is finite. In this case $X(a,b)$ 
(the ``excess'' number) is defined to be
the least $k$ such that $a\setminus b \subseteq k$.
Thus $a \setminus X(a,b) \subseteq b$, but if $X(a,b)>0$ then $X(a,b)-1 \in a
\setminus b$.

A pre-gap is a pair of
sequences $g=\{ (a_i\mid i \in I), (b_j \mid j \in J)\}$ where $I,\; J\subseteq
\omega_1$ are uncountable and $a_i,b_j\in [\omega ]^\omega$ are such that
\[ a_{i_0} \subseteq^* a_{i_1} \subseteq^* b_{j_1}\subseteq^* b_{j_0}\]
whenever $i_0 <i_1$  are in $I$ and $j_0 < j_1$ in $J$. 
In most cases
$I=J=\omega_1$. Given a pre-gap as above, 
and uncountable subsets $I'\subseteq I$ and $J'\subseteq J$, 
the restriction $g\restriction(I',J')$ of $g$ is the 
pre-gap $\{ (a_i\mid i \in I'), (b_j \mid j \in J')\}$. We write
$g\restriction I$ for $g\restriction(I,I)$.

An interpolation for a pre-gap $g$ is a set $x\subseteq \omega$ such that
$$a_i \subseteq^* x \subseteq^* a_j$$ for every $i$ and $j$. A pre-gap with
no interpolation is called a gap. A famous construction of
Hausdorff produces gaps in ZFC (which are now called Hausdorff gaps).
Specifically, a Hausdorff gap is a pre-gap $g =\{ (a_i\mid i \in \omega_1),
 (b_j \mid j \in \omega_1)\}$ such that
for every  $\alpha \in \omega_1$ and $n\in \omega$ the set
\[ \{ \beta \in \alpha \mid a_\beta \setminus n \subseteq b_\alpha \}\]
is finite.

 A {\em special} (or Kunen) gap is a pre-gap 
 $g =\{ (a_i\mid i \in \omega_1), (b_j \mid j \in \omega_1)\}$
such that for some $n_0\in \omega$:
\begin{enumerate}
\item $a_\alpha \setminus n_0 \subseteq b_\alpha$ for every $\alpha \in
\omega_1$, and
\item for all $\alpha < \beta  < \omega_1$,
\[ (a_\alpha \cup a_\beta)\setminus n_0 \not \subseteq b_\alpha \cap
b_\beta\]
(equivalently,  $a_\alpha\setminus n_0  \not \subseteq b_\beta$ or
$a_\beta\setminus n_0  \not \subseteq b_\alpha$).
\end{enumerate}

The interest in these definitions arises from
the fact
(not too difficult to prove) that these Hausdorff and Kunen 
pre-gaps  are gaps and remain gaps as long as $\omega_1$ is not
collapsed: they have no interpolation in
any extension in which $\omega_1$ remains uncountable.

In this paper we define two additional types of ``special''
gaps: $S$-Hausdorff gaps where $S\subseteq \omega_1$ is a stationary set, and
$C$-Hausdorff gaps, where $C$ is a ladder system over $S$.

The motivation for this work is the desire to find an example with
 gaps of the phenomenon in which $\omega_1$ is ``split'' in a certain 
 behavior on a stationary set $S\subset \omega_1$ and an opposite behavior on
  its complement $\omega_1 \setminus S$.

\begin{definition}
Let $S\subseteq \omega_1$ be a stationary set. A pre-gap
$g =\{ (a_i\mid i \in \omega_1), (b_j \mid j \in \omega_1)\}$
is $S$-Hausdorff iff for some closed
 unbounded (club) set $D \subseteq \omega_1$, for every $\delta \in S \cap D$
 and  for every sequence of ordinals
 $(i_n \in \delta \mid n \in \omega )$ increasing and cofinal in
 $\delta$
 \begin{equation}
 \label{Lim}
  \lim_{n\to \infty} X(a_{i_n},b_\delta ) = \infty.
  \end{equation}
 That is, for every $k$, there is only a finite number of $n\in \omega$ for 
which $a_{i_n} \setminus k \subseteq b_\delta$.
Since, for $\delta < \delta'$, $b_{\delta'} \subseteq^* b_\delta$,
 it follows that eventually 
$X(a_{i_n},b_\delta) \leq X(a_{i_n}, b_{\delta'})$, and hence 
(\ref{Lim}) holds for every $\delta' \geq \delta$ in $\omega_1$.
 If the pre-gap $g =\{ (a_i\mid i \in I), (b_j \mid j \in J)\}$ 
 is defined only on uncountable sets $I,J \subseteq \omega_1$, we
 can still define it to be $S$-Hausdorff if for some closed unbounded
 set $D\subseteq \omega_1$, for every $\delta\in S \cap D$, for every
 $j\in J \setminus \delta$, and for every increasing sequence of
 ordinals $ (i_n\in \delta \cap I \mid n \in \omega )$
 cofinal in $\delta$,
\begin{equation}
\label{Lim2}
 \lim_{n \to \infty} X(a_{i_n}, b_j)= \infty
 \end{equation}

 \end{definition}

 Clearly, every Hausdorff gap is an $\omega_1$-Hausdorff gap, and the
 closed unbounded set $D$ can be taken to be $\omega_1$.
 The converse of this also holds, in the sense that every
 $\omega_1$-Hausdorff gap contains a Hausdorff gap. For suppose that
 $g= \{ (a_i\mid i \in I), (b_j \mid j \in J)\}$ 
 is some  $\omega_1$-Hausdorff gap, and let $D\subseteq \omega_1$
 be the closed unbounded set given by the definition of $g$ as an 
  $\omega_1$-Hausdorff gap. Define $I' \subseteq I$ such that every
  two members of $I'$ contain a point from $D$ in between. We claim
  that $g' = \{ (a_i\mid i \in I'), (b_j \mid j \in J)\}$, 
 is a Hausdorff gap. Indeed, if $\alpha \in J$ then for every $n\in
  \omega$ the set $E= \{ \beta\in \alpha\cap I' \mid a_\beta \setminus n
  \subseteq b_\alpha \}$ is necessarily finite. For if not, then let
  $\delta$ be an accumulation point of $E$, 
  and let $\beta_i\in E$, for $i\in \omega$,
  be increasing and converging to $\delta$. Necessarily $\delta \in D$, 
  and  $a_{\beta_i} \setminus n \subseteq b_\alpha$ 
  shows that (\ref{Lim2}) does not hold.

\begin{prop}
\label{StatGap}
If $S\subseteq \omega_1$ is stationary, then
any $S$-Hausdorff pre-gap is a gap. (So that any pre-gap
containing an $S$-Hausdorff pre-gap is a gap.) 
\end{prop}
{\bf Proof.} Assume that this not so and let $x$ be an interpolation of
an
$S$-Hausdorff pre-gap $g= \{ (a_i\mid i \in I), (b_j \mid j \in J)\}$.
Then there is a fixed $n_0\in \omega$
 such that for unbounded  sets of indices $I'\subseteq I$
 $J'\subseteq J$,
for every $\alpha\in I',\; \beta \in J'$
\[ a_\alpha \setminus n_0 \subseteq x \setminus n_0 \subseteq b_\beta.\]
And as a consequence 
\begin{equation}
\label{E1}
 a_\alpha \setminus n_0 \subseteq b_\beta
\end{equation}
holds. Since $g$ is assumed to be $S$-Hausdorff there exists a
closed unbounded set $D\subseteq \omega_1$ as in the definition. We
may assume that every $\delta\in D$ is an accumulation point of $I'$.
Take now a limit $\delta \in S\cap D$ 
  and any sequence  $i_n\in I'\cap \delta$
increasing to $\delta$. Take any $j\in J' \setminus \delta$.
 Then equation (\ref{E1})
implies that 
$X(a_{i_n},b_\delta)\leq n_0$, which is a contradiction.\qed

A stronger notion than that of being $S$-Hausdorff can be defined if
the rate at which the sequences in (\ref{Lim}) tend to infinity
is uniform.
For this we must recall the definition of a scale or ladder system on
a stationary set.

If $S\subseteq \omega_1$ is a stationary set, then a ladder (system) over $S$
is a sequence $C =  ( c_\alpha \mid \alpha \in S\ \mbox{\it is a limit ordinal} )$
such that every $c_\alpha=( c_\alpha(n)\mid n \in \omega )$ 
is an increasing, cofinal in $\alpha$
$\omega$-sequence. 
\begin{definition}
\label{DHG}
For a ladder system $C$ over $S$, we say that a pre-gap
$g = \IJpregap$ is $C$-Hausdorff iff for some closed-unbounded set
$D\subseteq \omega_1$ for all $\delta \in S\cap D$ and $j\in J
\setminus \delta$ 
there is $k\in \omega$ such that for
every $n\geq k$ in $\omega$, if $i \in I\cap (\delta \setminus c_\delta(n))$, 
then $X(a_i, b_j)> n$.
\end{definition}

Every $C$-Hausdorff gap (where $C$ is a ladder system over a
stationary set $S$) is $S$-Hausdorff. Our aim in this paper is to prove the
following consistency result.
\begin{theorem}
Assume G.C.H for simplicity.
Suppose that $\kappa$ is a cardinal such that $\cf(\kappa)> \aleph_1$. 
Let $S$ be a stationary co-stationary subset of $\omega_1$.
Then there is a c.c.c  poset of size 
$\kappa$ such that in every generic extension made via $P$ 
$2^{\aleph_0}=\kappa$ and
the following hold.
\begin{enumerate}
\item For every ladder system $C$ over $S$, every 
gap contains a subgap that is $C$-Hausdorff.
\item For every ladder system $E$ over $\omega_1\setminus S$ there
is a gap $g$ with no subgap that is $E$-Hausdorff.
\end{enumerate}
\end{theorem}

\section{Gaps introduced by forcing}
\label{GF}
Gaps can be created by forcing with finite conditions (a method due to
Hechler \cite{Hechler}). These gaps are not $S$-Hausdorff for any 
stationary set, as we are going to see.

If $f\in 2^n$ ($f$ is a function defined on $n$ 
with range included in $\{ 0,1 \}$) then $f$ is a
characteristic function and we let $[ f ] = \{ k \mid f(k) =1 \}$ be
the subset of $n$ represented by $f$.

Let $(I,<_I)$ be any ordering isomorphic to $\omega_1 + \omega_1^\ast$. 
For example take $I = (\omega_1 \times \{ 0 \}) \cup (\omega_1 \times\{ 1 \})$
with $\langle \alpha,0\rangle <_I \langle \beta,0\rangle
<_I \langle \beta,1 \rangle <_I \langle \alpha, 1 \rangle$ 
whenever $\alpha < \beta < \omega_1$.

Define the poset $P$ by $p \in P$ iff $p$ is a finite function defined on
 $I$ and
such that:
\begin{enumerate}
\item  For some $n$ (called the ``height'' of $p$) $p(i) \in 2^n$ for every $i\in
\dom(p)$. (The height of the empty function is defined to be $0$.)
\item For every $\alpha \in \omega_1$, $\langle \alpha, 0 \rangle \in\dom(p)$
 iff $\langle \alpha, 1 \rangle \in\dom(p)$, and in this case $[
p(\langle \alpha, 0 \rangle) ] \subseteq [ p(\langle \alpha, 1
\rangle) ]$.
\end{enumerate}
  The intuition behind this definition is that for 
$\alpha \in \omega_1$, $p(\langle \alpha,0\rangle)$ 
will ``grow'' to become $a_\alpha$, and $p(\langle \alpha,
1\rangle)$ will finally become $b_\alpha$, as $p$ runs over the generic filter. 
So that $(\langle a_\alpha \mid
\alpha \in \omega_1\rangle, \langle b_\alpha\mid \alpha\in \omega_1\rangle)$
will be the generic gap with the additional property that $a_\alpha
\subseteq b_\alpha$ for every $\alpha$.
 The ordering of $P$ reflects this intuition as
follows.

For $p_1,p_2\in P$ define $p_1 \leq p_2$ ($p_2$ extends $p_1$) iff 
\begin{enumerate}
\item $d_1 = \dom(p_1) \subseteq d_2=\dom(p_2)$, and for every 
$i\in d_1$,
$p_1(i) \subseteq p_2(i)$ (so $\height(p_1)\leq \height(p_2)$).
\item For every $i,j\in \dom(p_1)$, if $i <_I j$ then
\[ [p_2(i)] \setminus [ p_1(i)] \subseteq [ p_2(j)].\]
\end{enumerate}

It is easy to see that any condition in $P$ has extensions
 with arbitrarily large
height and with domains that extend arbitrarily over $I$.
In fact, given $i\in \dom(p)$ and $k\in \omega$ above height $p$, we can require
that the extension $p'$ puts $k$ in $[p'(i)]$.

If $\alpha \in \omega_1$, we can write $\alpha \in \dom(p)$ instead of
$\langle \alpha, 0\rangle\in \dom(p)$ (which is equivalent to 
$\langle \alpha, 1\rangle\in \dom(p)$). So $\dom(p)$ has two meanings, and 
the context decides if it means a set of ordinals or a set of pairs.

Suppose that $A \subseteq I$ is such that $\langle \alpha, 0 \rangle\in A$ iff
$\langle \alpha, 1 \rangle\in A$. Let $P_A$ be the subposet of $P$
consisting of all conditions $p$ such that $\dom(p) \subseteq A$.
If $p\in P$ then $p\restriction A \in P_A$ and $p\restriction A \leq p$.
We prove some additional properties of this restriction map taking $p$ to
 $p\restriction A$.

In the definition of $p\leq q$ what really counts is the
restriction of $q$ to the domain of $p$. That is,
 $p \leq q $ iff $p \leq q \restriction \dom(p)$.
It follows that $p \leq q$ implies that $p \restriction A \leq q
\restriction A$. It also follows that if $p$ and $q$ are conditions 
such that for $C =\dom(p) \cap \dom(q)$, 
$p\restriction C = q \restriction C$, then $p$ and $q$ are compatible.
In fact, in this case, $p\cup q$ is the minimal extension of $p$ and $q$.

Suppose that $\dom(p) = \dom(q)$. Then $p$ and $q$ are compatible
in $P$ iff $p \leq q$ or $q \leq p$.

For compatible conditions $p$ and $q$, we define a canonical extension
$p\vee q$ of both $p$ and $q$. However, $P$ is not a lattice and 
$p\vee q$ is not the minimum of all extensions of $p$ and $q$. 
To define it, we first 
make an observation. Consider $C = \dom(p) \cap \dom(q)$. 
Then $p\restriction C$ and $q\restriction C$ are comparable in $P_C$ 
(since they are compatible and have the same domain), 
and hence we can assume without loss of generality that 
 $q\geq p\restriction A$ and 
$n=\height(q)\geq m = \height(p)$ where $A = \dom(q)$
 (the restriction on the heights is 
needed only in case $p\restriction A = \emptyset$ since it follows
from  $q\geq p\restriction A$ otherwise). Then $r=p \vee q$ is
defined as follows on $\dom(p) \cup \dom(q)$,
 and it will be evident that $p \vee q$ is an extension of
$p$ and $q$.

For $i\in \dom(q)$ define $r(i) = q(i)$. For $i\in \dom(p) \setminus A$ define
$r(i)\in 2^n$ by the following two conditions:
\begin{equation}
 p(i) \subseteq r(i).
 \end{equation}
\begin{equation}
\label{cc}
[ r(i)] \setminus [ p(i) ] = \bigcup \{ [q(k)]\setminus m \mid k <_I i\ \mbox{and } k
\in A \cap \dom(p) \}.
\end{equation}
This definition makes sense since $A\cap \dom(p) \subseteq \dom(q)$.

It is clear that $r\in P$, $\dom(r)=\dom(p) \cup \dom(q)$ and
$r\restriction A = q$. We prove that $r \geq p$.
Clause 1 in the definition of extension is obvious, and we have to
check clause 2.
 Suppose that $i,j \in
\dom(p)$ and $i <_I j$. We have to show that
\begin{equation}
\label{aim}
 [r(i)] \setminus [p(i)] \subseteq [r(j)].
\end{equation}
So consider any $a \in [r(i)] \setminus [p(i)]$.

\noindent
{\bf Case 1:} $i\in A$. Then $r(i) = q(i)$. 
If $j\in A$ as well, then (\ref{aim}) follows from
our assumption that $q\geq p\restriction A$, and since $r(i) = q(i)$, $r(j)
= q(j)$ in this case. If, on the other hand, $j \not \in A$, then
\[ [r(j) ] \setminus [p(j)] = \bigcup \{ [q(k)]\setminus m \mid k <_I j\ \mbox{ and } k
\in A \cap \dom(p)\}\]
by the definition of $r$.
Since $i\in A \cap \dom(p)$, $i<_I j$, and $a\in q(i) \setminus m$, 
$a \in [ r(j) ] \setminus [ p(j) ]$ as required.

\noindent
{\bf Case 2:} $i \not \in A$. Then $i \in \dom(p) \setminus A$ and
(\ref{cc}) implies that for some $k\in A\cap \dom(p)$ such that $k <_I i$, $a \in
[ q(k) ]\setminus m$. 
Then $k <_I j$,  both indices are in $\dom(p)$, and $k \in A$, which brings
us back to Case 1. \qed

This argument has the following corollary.
\begin{cor}
Suppose that $p_1,p_2\in P$ and $C=\dom(p_1) \cap \dom(p_2)$ are
such that $p_1\restriction C \geq p_2\restriction C$ and 
$\height(p_1)\geq  \height(p_2)$. Then $p_1
\vee p_2$ can be formed (an extension of $p_1$ and $p_2$).
\end{cor}
{\bf Proof.} Define $A = \dom(p_1)$. Then $p_1 \geq p_2\restriction A$ (because
$p_1 \geq p_1\restriction C \geq p_2\restriction C = p_2\restriction A$).
 So $r = p_1 \vee p_2$ can be formed.\qed

\begin{lemma}
$P$ satisfies the c.c.c. In fact if $\{ p_\alpha \mid \alpha\in S\} \subseteq
P$ where $S\subseteq \omega_1$ is stationary, then for some stationary set
$S'\subseteq S$, every finite set of conditions in $\{ p_\alpha \mid \alpha
\in S'\}$ is compatible. (This is Talayaco's condition \cite{Talayaco}.)
\end{lemma}
{\bf Proof.} If $p,q\in P$ have the same height 
and for $C=\dom(p) \cap \dom(q)$ it happens that
$p\restriction C= q\restriction C$, 
then $p\cup q$ is an extension of $p$ and $q$. Hence a
$\Delta$-system argument works here.\qed

If $G\subset P$ is some generic filter over $P$, define for
every $\alpha \in \omega_1$ $a_\alpha = \bigcup \{ [ p(\langle\alpha,0
\rangle)] \mid p \in G \}$, and $b_\alpha = \bigcup \{ [
p(\langle\alpha,1\rangle)] \mid p \in G \}$. A standard density argument
shows that $g$ is a pre-gap, and we denote it as $g$.
\begin{lemma}
\label{Gap}
The generic pre-gap $g$ is a gap.
\end{lemma}
{\bf Proof.} Suppose that $x\in V^P$ is a name, forced to be an
interpolation for the generic pre-gap $g$. For every $\alpha \in \omega_1$ 
find a condition $p_\alpha\in P$ and a number $n_\alpha\in \omega$ such
that
\begin{equation}
\label{PN}
p_\alpha \forcesP a_\alpha \setminus n_\alpha \subseteq x\setminus n_\alpha
 \subseteq b_\alpha.
\end{equation}
Then for some stationary set $S\subseteq \omega_1$, and some fixed $n\in
\omega$, $n=n_\alpha$ for every $\alpha \in S$, and the sets
$\dom(p_\alpha)$ form a $\Delta$-system  with core $C\subset I$ (a finite set). We
also assume that $p_\alpha \restriction C$ is fixed for $\alpha \in
S$. For $\alpha < \beta$, both in $S$ and above the ordinals involved
in $C$, 
consider $p_\alpha$ and $ p_\beta$.
 Pick any $k\geq n$ such that $k\geq  \height(p_\alpha)$ as well.
Let $i = \langle \alpha , 0 \rangle$, and $j=\langle \beta, 1
\rangle$. We shall find an extension $r$ of $p_\alpha$ and $p_\beta$
such that $r(i)(k) = 1$ and $r(j)(k) = 0$. Then $r\forces k \in
a_\alpha \wedge k \not \in b_\beta$. But this contradicts (\ref{PN}).

To define $r$, define first an extension 
$p'_\alpha \geq p_\alpha$ by requiring
that $p'_\alpha(i)(k) = 1$ and $[p'_\alpha(\langle \gamma,0\rangle)] =
[p_\alpha(\langle \gamma,0\rangle)]$ for every $\langle \gamma,0\rangle \in C$.
This is possible since $i$ is never $<_I$ below $\langle \gamma, 0 \rangle\in C$.
Now $p'_\alpha$ extends $p_\beta\restriction C$ and hence $r = p'_\alpha \vee
p_\beta$ can be formed. Since the only members of $C$ below $j$ (in $<_I$) are of
the form $\langle \gamma,0\rangle$, it follows that $[r(j)]=[p_\beta(j)]$. Thus
$r(j)(k)=0$.
\qed

The following lemma implies that 
if $G$ is a $(V,P)$-generic filter, $g$ the generic gap,
 and $U\in V[G]$ is any stationary subset of $\omega_1$ in the
extension, then no uncountable restriction of 
$g$ is $U$-Hausdorff.

\begin{lemma} 
\label{Prop}
The following holds in $V^P$ for the generic gap $g=
\{ (a_i \mid i \in \omega_1), (b_j\mid j \in \omega_1) \}$. 
If $J, K \subseteq \omega_1$ are unbounded, then
there is  a club set $D_0\subseteq \omega_1$ such that for every 
$\delta \in D_0 $ and
$k\in K\setminus \delta$ there are $m\in \omega$ and a sequence
$j(n)\in \delta \cap J$ increasing and cofinal in $\delta$ such that
 $a_{j(n)}\setminus m \subset b_k$ for all $n\in \omega$.
\end{lemma}
{\bf Prof.} Let $J, K \in V^P$ be names forced by every condition in $P$ to be
unbounded  subsets of $\omega_1$. Define in $V^P$ the following set
$D_0\subseteq \omega_1$:  $\delta \in D_0$ if and only if
$\delta\in \omega_1$ is a limit ordinal such that:

\[
\begin{minipage}[t]{110mm} 
 \it for all $k \in K \setminus \delta$
{ there is some $m\in \omega$ and an increasing, cofinal in } 
$\delta$  { sequence }
$j(n)\in \delta \cap J$ \mbox{ with } $a_{j(n)}\setminus m \subseteq b_k$.
\end{minipage}
\]

We want to prove that $D_0$ contains a closed unbounded subset of 
$\omega_1$, and assume that it does not. So $R=\omega_1 \setminus D_0$ is
(forced by some condition to be) 
stationary in  $V^P$, and hence the set, defined in $V$, of ordinals that
are potentially in $R$ is stationary in $V$.
Namely, the set $R_0\subset \omega_1$ of ordinals forced by some
condition  to be in $R$ is stationary. For every $\delta\in
R_0$ pick a condition $p_\delta$ that forces $\delta \not \in D_0$.
By extending $p_\delta$ we can find
some $k_\delta\geq \delta$ such that 
$$p_\delta\forcesP k_\delta \in K\ \mbox{shows that } \delta \not
\in D_0.$$
By extending $p_\delta$ again, we can find some $j_\delta\in
\omega_1 \setminus \delta$ forced by $p_\delta$ to be in $J$ (which is possible
since $J$ is supposed to be unbounded in $\omega_1$). 
If necessary, a further extension ensures that both $j_\delta$ and $k_\delta$ are in the domain of $p_\delta$.
Now there exists some $m = m_\delta \in \omega$ such that
$p_\delta \forces\; a_{j_\delta} \setminus m \subseteq b_{k_\delta}$
(the height of $p_\delta$ will do). We can extend $p_\delta$ once
again and find $f(\delta) < \delta$ such that 
\begin{equation}
\label{FEq}
p_\delta \forcesP \mbox{ there is no } j \in J,\ f(\delta)< j <
\delta,\ \mbox{for which } a_j \setminus m \subseteq b_{k_\delta}.
\end{equation}

We may assume that, for a stationary set
$T\subseteq R_0$, the domains of $p_\alpha$, for $\alpha \in T$, form a
$\Delta$ system, that  they all have the same
height, say $n$, and the same restriction to the core. 
  We also assume that the functions 
$p_\alpha(\langle j_\alpha, 0 \rangle): n \to \{ 0 , 1 \}$
 do not depend on $\alpha$, and
  that $f(\alpha)$ and $m= m_\alpha$ are fixed on $T$ ($m\leq n$).
Now by Talayaco's
 chain condition for $P$, there is a stationary $T'\subseteq T$ such
that for every $\alpha,\beta\in T'$, $p_\alpha \vee p_\beta$ is a 
common extension. Pick some
$\alpha \in T'$ that is an accumulation point of $T'$ (and such that for every
$\beta < \alpha$, $j_\beta < \alpha$). Then find
$\beta<\alpha$, $\beta\in T'$ such that $f(\alpha) < \beta$.
Then (as we shall see)
\[ p_\beta\vee p_\alpha \forces a_{j_\beta} \subseteq a_{j_\alpha},\ 
j_\beta\in J,\
\mbox{\it and } \alpha > j_\beta > f(\alpha).\]
Yet
\[ p_\alpha \forces\; a_{j_\alpha} \setminus m \subseteq b_{k_\alpha},\]
and this is a contradiction to (\ref{FEq}).
Why does $p_\beta \vee p_\alpha$ force $a_{j_\beta} \subset a_{j_\alpha}$?
 Because the functions $p_\alpha(\langle j_\alpha,0 \rangle)$ and
  $p_\beta(\langle j_\beta,0 \rangle)$ are the same, they describe
  $a_{j_\alpha}\cap n$ and $a_{j_\beta}\cap n$, so $p_\beta\vee p_\alpha$ forces
  $a_{j_\beta}\subset a_{j_\alpha}$.
\qed

\section{Specializing  pre-gaps on a ladder system}
\begin{theorem}
\label{Sp}
For every ladder system $C$ over a stationary set $S\subseteq \omega_1$,
and gap $g$, there is a c.c.c forcing notion $Q = Q_{g,C}$ such that in $V^Q$ a
restriction of $g$ to some uncountable set is $C$-Hausdorff 
(and hence $S$-Hausdorff).  In fact $Q$ 
satisfies a stronger property than c.c.c, the polarized chain condition, 
which we shall define later.
\end{theorem}
{\bf Proof.}
Fix for the proof a ladder system 
$C = \langle c_\delta \mid \delta \in S \rangle$ over a stationary set $S
\subseteq \omega_1$ consisting of limit ordinals,
 and a pre-gap $g=\pregap$.
The forcing poset $Q =Q_{g,C}$ defined below is designed to make an uncountable restriction of 
$g$ into a $C$-Hausdorff gap.

Define $p\in Q$ iff $p=(w,s)$ where
\begin{enumerate}
\item $w\in [ \omega_1 ] ^{<\aleph_0}$ (i.e. a finite subset of 
$\omega_1$), and 
\item $s\in [ S ] ^{<\aleph_0}$.
\end{enumerate}
If $p\in Q$ then we write $p=(w^p,s^p)$ for the two components of $p$.

The ordering $p\leq q$ ($q$ extends $p$) is defined by
\begin{description}
\item[a.]  $w^p \subseteq w^q$, $s^p\subseteq s^q$, and
\item[b.] If $\delta \in s^p$ and $i\in w^p$ are such that $\delta \leq i$,
then for every $j\in (w^q\setminus w^p) \cap \delta$,

\[ a_j \setminus \mid c_\delta \cap j\mid\;  \not \subseteq b_i.\]
\end{description}
Or, equivalently, $X(a_j,b_i) > |c_\delta \cap j |$.
It is easy to check that this is indeed an ordering defined on $Q$.

If $G$ is generic over $Q$, define 
$W = \bigcup \{ w \mid \exists s (w,s)\in G \}$. 
We will prove that if $g$ is a gap then $Q$ satisfy the c.c.c. 
 So $\omega_1$ is preserved.
Clearly, if $p=(w,s)$ is a condition, then for any $\sigma\in S$,
$(w,s\cup\{ \sigma\})$ extends $p$, and if $j\in \omega_1$ and $j>\max(w)$, then
$(w\cup\{ j \},s)$ extends $p$.
(If, however, $j < \max(w)$, then $(w\cup \{ j \}, s)$ may be
incompatible with $(w,s)$.) It follows that $W$ is unbounded in
$\omega_1$ and $\{ ( a_i \mid i \in W), (b_j \mid j \in
W)\}$ is $C$-Hausdorff.

So the generic filter over $Q$ selects an
unbounded in $\omega_1$ restriction of $g$ that is $C$-Hausdorff.

If $p=(w,s)$ is a condition then for every $\alpha\in \omega_1$ 
the restriction $p\restriction \alpha = ( w\cap \alpha,
s\cap \alpha)$ is defined. Clearly $p\restriction \alpha \leq p$.

If $p=(w,s)$ and $q=(v,r)$ are conditions 
in $Q$ then define $p\cup q=(w\cup v, s\cup r)$.
If $p$ and $q$ are compatible in $Q$, then $p\cup q\in Q$ is the least upper
bound of $p$ and $q$.

The following lemma describes a situation in which the compatibility of 
$p_1$ and $p_2$ can be deduced. 
This is the situation resulting when $p_1$ and $p_2$ come from a
$\Delta$-system, with core fixed below $\gamma$, and such that $p_1$ is bounded by
some $\alpha$ such that 
the domain of $p_2$ has empty intersection with the ordinal interval
$[\gamma,\alpha]$. The proof is straightforwards.
\begin{lemma}
\label{l32}
Suppose that
\begin{enumerate}
\item $p_1=(w_1,s_1)$ and $p_2=(w_2,s_2)$ are in $P$.
\item 
$\gamma < \alpha < \omega_1$ are such that 
\begin{enumerate}
\item
$ w_1 \subseteq \alpha,$
and $p_1\restriction \gamma$ is compatible with $p_2$.
\item
$ w_2 \cap \alpha \subset \gamma$, and $s_2 \cap (\alpha + 1) \subset \gamma$.
$p_2\restriction \alpha = p_2 \restriction \gamma$ is compatible with $p_1$.
\end{enumerate}

\item Define
\[ A = \bigcap \{ a_i \mid i \in w_1 \setminus \gamma \}\]
\[ B = \bigcup \{ b_j \mid j \in w_2 \setminus \gamma \}\]
and suppose that there is $n \in A \setminus B$ such that, for every
$\delta \in s_2 \setminus \alpha$, $n > | c_\delta \cap
\alpha|$.
\end{enumerate}
Then  $p_1$ and $p_2$ are compatible.
\end{lemma}
{\bf Proof.} Form $p = p_1 \cup p_2$ and prove that $p_1, p_2 \leq p$. $p_1\leq
p$ is immediate. As for $p_2 \leq p$, observe that $X(a_i,b_j) > | c_\delta
\cap \alpha |$ for every $i \in w _1 \setminus \gamma$, $j \in w_2\setminus
\gamma$, and $\delta \in s_2 \setminus \gamma$. \qed

The following simple lemma is used in proving that $Q$ is a c.c.c
poset.
\begin{lemma}
\label{simp}
Suppose $g=\{ (A_i\mid i \in I), ( B_j\mid j \in J) \}$ is a pre-gap
such that for every $i\in I$ and $j\in J$, $i<j$ implies that $A_i\subseteq
B_j$. Then $g$ is not a gap.
\end{lemma}
{\bf Proof.}
By throwing away a countable set of indices from $J$ we can assume for
every $n\in \omega$ that if $n\not\in B_j$ for some $j$, then $n\not \in
B_j$ for uncountably many $j$'s. Define then $x=\bigcup_{i\in I} A_i$. Then
$x\subseteq B_j$ for every $j$, because otherwise there are some $i\in I$,
$j\in J$, and $n\in \omega$ such that $n\in A_i\setminus B_j$. But then we
may find uncountably many indices $j'$ such that 
 $n\not\in B_{j'}$ and in particular there is such $j' > i$. 
 Thence $A_i\not\subseteq B_{j'}$,
contradicting our assumption. \qed

\begin{theorem}
\label{PCC}
Suppose that the domain of our ladder system 
$C$, namely $S$, is co-stationary. 
\begin{enumerate}
\item
If $g$ is a gap then $Q=Q_{g,C}$ satisfies the c.c.c.
\item 
Suppose that 
$T_1$, $T_2 \subseteq \omega_1 \setminus S$ are stationary sets  
and $\ovp = ( p^\ell_\delta \mid \delta \in T_\ell)$, for $\ell =
1,2$ , are
two sequences of conditions in $Q$ such that, for some fixed $p^*\in
Q$, $p^* \geq p_\delta^1\restriction \delta ,\  p^2_\mu\restriction \mu$, for every
$\delta \in T_1$ and $\mu \in T_2$, is such that $p^*$ is compatible with
every $p_\delta^1$ and with every $p^2_\mu\restriction \mu$. 
Then there are stationary subsets
$T_1'\subseteq T_1$ and $T_2'\subseteq T_2$ such that, for every
$\alpha_1\in T'_1$ and $\alpha_2\in T'_2$, if $\alpha_1 < \alpha_2$
then $p^1_{\alpha_1}$ and $p^2_{\alpha_2}$ are compatible in $Q$.
\end{enumerate}
\end{theorem}
{\bf Proof.}
We prove {\it 2} since the proof of {\it 1}
is similar. For any condition $p=(w,s)$
define $\dom(p)=w\cup s$.  
Suppose that $\dom(p^*)\subseteq \gamma$. Then 
$\dom(p_\delta^1) \cap \delta \subseteq \gamma$, and $\dom(p_\mu^2)\cap \mu
\subseteq \gamma$,
for every $\delta \in T_1$ and $\mu \in T_2$. 
We may assume that if $i<j$ then $\dom(p^\ell_i) \subset \cap\;
(\dom(p^m_j)\setminus \gamma)$ for $\ell,m\in\{1,2\}$.

Since $\delta \in T_\ell$ implies that $\delta \not\in S$, it follows 
for $p_\delta^\ell = (w,s)$ that
the ladder sequence $c_i$ for any $i \in s$ is bounded below
$\delta$. So the finite union 
\[ \delta\cap \bigcup\{ c_i\mid i \in s^{p^\ell_\delta} \setminus \delta \}\]
is bounded below $\delta$. Using Fodor's lemma we may even assume that this
intersection is bounded below $\gamma$ (extend $\gamma$ if necessary) and has a
fixed finite cardinality.

For every $\delta \in T_1$ define
\[ A_\delta = \bigcap \{ a_i \mid i\in w^{ p^1_\delta}\setminus
\gamma\}\]
Similarly, for $\delta \in T_2$ define
\[ B_\delta = \bigcup\{ b_i \mid i\in w^{ p^2_\delta}\setminus
\gamma\}. \]
Clearly, any interpolation for 
$G=\{ (A_\delta\mid \delta\in T_1\},\{ B_\delta\mid \delta\in
T_2\})$ is also an interpolation for $g$, and hence $G$ is a gap.

Let $k\in \omega$ be such that for every $\delta\in T_\ell$, if
$p^\ell_\delta=(w,s)$ and $\alpha\in s\setminus \gamma$, then $\mid
c_\alpha\cap \delta\mid < k$.

Now we find a stationary set $T'_1\subset T_1$ such that for every $n\in
\omega$ if $n\in A_\delta$ for some $\delta \in T'_1$ then $n\in
A_\delta$ for a stationary set of $\delta$'s in $T'_1$. 
Simply throw away countably
many non-stationary sets from $T_1$. Similarly, find  a stationary
$T'_2\subseteq T_2$ such that if $n\not \in B_\delta$
for some $\delta\in T'_2$ then $n\not \in B_\delta$ for a stationary set
of $\delta\in T'_2$. 

Now lemma \ref{simp} gives $\alpha_1\in T'_1$ and $\alpha_2\in T'_2$
with $\alpha_1<\alpha_2$ such that $A_{\alpha_1}\setminus k \not \subseteq 
B_{\alpha_2}$. If we pick $n\in A_{\alpha_1}\setminus B_{\alpha_2}$ such
that $n\geq k$ then there are stationary sets $T^{''}_1\subseteq T'_1$ and
$T^{''}_2\subseteq T'_2$ such that  $n\in A_{\alpha_1}\setminus B_{\alpha_2}$
{\em for every} $\alpha_1\in T^{''}_1$ and $\alpha_2\in T^{''}_2$.
Hence if $\alpha_1\in T^{''}_1$, $\alpha_2\in T^{''}_2$, and
$\alpha_1<\alpha_2$, then $p^1_{\alpha_1}$ and $p^2_{\alpha_2}$ are
compatible in $Q$ by lemma \ref{l32}.\qed

\subsection{Polarized chain condition}
Theorem \ref{PCC} shows that the poset $Q_{g,C}$ for a gap $g$ and ladder $C$
over a stationary co-stationary set $S$ satisfies some
kind of a chain condition, suitable for two sequences indexed by stationary
subsets of $\omega_1\setminus S$. We formulate this condition in general
and later  prove 
that the iteration with finite support  preserves this condition.

\begin{definition}
Let $T\subseteq \omega_1$ be a stationary set. 
A c.c.c poset $P$ satisfies the polarized chain condition ({\em
p.c.c})
for $T$ if it satisfies the following
requirement. Suppose that
\begin{enumerate} 
\item
$$\ovp^\ell=(p^\ell_\delta\mid \delta\in T_\ell)\ \mbox{for } \ell
= 1,2$$
are two sequences of conditions in $P$, where $T_\ell\subseteq T$
are stationary for $\ell = 1,2$.
\item $p^*\in P$ is such that for each $\ell = 1,2$
\[ p^*\forcesP \{ \delta\in T_\ell\mid p^\ell_\delta\in G_P\}\ \mbox{is
stationary in} \ \omega_1,\]
where $G_P$ is the name of the generic filter over $P$.
\end{enumerate}
Then there are stationary sets $T'_\ell\subseteq T_\ell$ for $\ell=1,2$
such that $p^1_{\alpha_1}$ and  $p^2_{\alpha_2}$ are compatible in $P$
whenever $\alpha_1< \alpha_2$ are in $T'_1$ and $T'_2$ respectively.
\end{definition}

We want to prove that if $g$ is a gap and $C$ a ladder over a stationary
set $S$ such that $T=\omega_1 \setminus S$ is also stationary, then
$Q=Q_{g,C}$ satisfies the p.c.c. for $T$.
The problem is that if $p^*$ is as in the p.c.c. definition then it is not
necessarily of the form to which theorem  \ref{PCC} is applicable, and so
we need some argument to deduce that $Q$ is p.c.c.


Recall that every club (closed unbounded) subset of $\omega_1$
in a generic extension of $V$ made via a c.c.c poset
contains a club subset in $V$. The following property of
c.c.c posets is also needed.
\begin{lemma} 
\label{CCC}
Let $P$ be a c.c.c poset. Suppose that $T\subseteq \omega_1$
is stationary, and $\langle p_\alpha \mid \alpha \in T
\rangle$ is a sequence of conditions in $P$ indexed along
$T$. Then there exists some $p_\alpha$ such that
\[ p_\alpha \forces \{ \beta \in T \mid p_\beta \in G \}\ \mbox{is
stationary}.\]
In fact, the set of these $\alpha$'s is stationary in $\omega_1$.
\end{lemma}
{\bf Prof.} Assume that this is not the case and, for some club
$D\subseteq \omega_1$, for every
$\alpha \in T\cap D$ there is a club set $C_\alpha$ (necessarily in
$V$) and an extension $p_\alpha'$ of $p_\alpha$ such that 
\begin{equation}
\label{Cl}
p_\alpha' \forces \beta \in C_\alpha \cap T \rightarrow p_\beta \not
\in G.
\end{equation}

Let $C=\{ \beta\in \omega_1\mid (\forall \alpha < \beta)
\beta \in C_\alpha\}$ be the diagonal intersection of these
club sets. Then $C$ is closed unbounded in $\omega_1$. Take a maximal
antichain (surely countable) from the set of extensions 
$\{ p_\alpha' \mid \alpha\in T\cap D \}$,
and let $\alpha_0$ be an index in $T\cap C \cap D$ 
higher than all indexes of
this countable antichain. Then $p'_{\alpha_0}$ is compatible with some
$p'_\alpha$ with $\alpha < \alpha_0$. 
But $\alpha_0\in C_\alpha$ which leads to a contradiction
since $p'_{\alpha_0}$ forces that $p_{\alpha_0}\in G$, 
 and $p'_\alpha$ forces that $p_{\alpha_0}\not \in G$
(by \ref{Cl}).\qed

Now we prove that $Q$ is p.c.c. for $T = \omega_1 \setminus S$.
\begin{lemma}
\label{LCC}
If $C$ is a ladder system over $S$, and $T=\omega_1 \setminus S$
is stationary, then, for any gap $g$, $Q_{g,c}$ is p.c.c. over $T$.
\end{lemma}
{\bf Proof.}
Suppose that
$T_1,T_2\subseteq T$ are stationary, and $\ovp^1$, $\ovp^2$ are two sequences
of conditions indexed along $T_1$ and $T_2$. Let $q^*\in Q$ be such that for
$\ell = 1 ,2$
\begin{equation}
\label{EQ}
q^* \forcesQ \{ \delta \in T_\ell \mid p_\delta^\ell \in G \}\ \mbox{is
stationary in }\ \omega_1.
\end{equation}

We claim first that we may assume that $q^* \leq p_\delta^1$ for every
$\delta\in T_1$. This can be achieved as follows. First, observe that the set
of $\delta\in T_1$ for which $p_\delta^1$ and $q^*$ are compatible is
stationary. Then use Fodor's theorem on this stationary set to fix
$p_\delta^1\restriction \delta$. Rename $T_1$ to the the resulting stationary
set. Redefine $p_\delta^1$ as $p_\delta^1\cup q^*$, and finally apply lemma
 \ref{CCC} to obtain $\delta_0$ such that
 \[ p_{\delta_0}^1 \forcesQ \{ \delta \in T_1 \mid p_\delta^1\in G \} \
 \mbox{is stationary in }\ \omega_1.\]

 Now $q^*$ is as required. Since $q^*$ extends the original  $q^*$, it still
 satisfies (\ref{EQ}) with respect to $T_2$. Repeat this procedure for $T_2$,
 and obtain two sequences to which Theorem \ref{PCC} is applicable.
 \qed

We note here for a possible future use 
 a stronger form of polarized chain condition
(strong-p.c.c) which is not used in this paper.
\begin{definition}
Let $T\subseteq \omega_1$ be stationary. A poset $P$ is said to satisfy the
 strong-p.c.c over $T$ 
if whenever two sequences are given
$$\ovp^\ell=(p^\ell_\delta\mid \delta\in T_\ell)\ \mbox{for } \ell
= 1,2$$ 
of conditions in $P$, where $T_\ell\subseteq T$
are stationary for $\ell = 1,2$, and for some $p^*\in P$, for every  $\ell = 1,2$,
\[ p^*\forcesP \{ \delta\in T_\ell\mid p^\ell_\delta\in G_P\}\ \mbox{is
stationary in} \ \omega_1,\]   
then there are stationary subsets $T_\ell'\subseteq T_\ell$ for $\ell = 1,2$, and
conditions $q_\delta \geq p^2_\delta$ for $\delta \in T'_2$ such that:

\begin{minipage}[t]{110mm}
For every $\delta\in T_2'$ and $q\in P$ such that $q_\delta\leq q$ there exists $\alpha <
\delta$ such that for every $\beta$ that satisfies $\alpha < \beta\in T_1'\cap \delta$
\[ q \ \mbox{and}\ p_\beta^1\ \mbox{are compatible in}\ P.\]
\end{minipage}
\end{definition}

\section{Iteration of p.c.c posets}

Our aim in this section is to prove that the iteration with finite support of
p.c.c  posets is again p.c.c. 
It is well known (by Martin and Solovay \cite{MartinSolovay}) that since each
of the iterands satisfies the countable chain condition the
iteration is again c.c.c, but we have to prove the preservation
of the polarized property.

Our posets are separative (and if not, they can be made separative). A poset
is separative iff $p \not \leq q$ implies that some extension of $q$ is 
incompatible with $p$.

\begin{lemma} 
\label{Pres}
Suppose that $P$ is a p.c.c. poset, and that $Q\in V^P$ is
(forced by every condition in $P$ to be) a p.c.c. poset. Then
the iteration $P\ast Q$ satisfies the polarized chain condition
too.
\end{lemma}
{\bf Proof.} Suppose that $(p^\ell_\delta, q^\ell_\delta)\in
P\ast Q$  are given for $\delta \in T_\ell \subseteq T$ and for
$\ell = 1,2$, such that for some condition $(p,q)\in P\ast Q$
 \begin{equation}
 \label{St} 
 (p,q) \forces \{ \delta \in T_\ell \mid (p^\ell_\delta, q^\ell_\delta) \in
G_{P\ast Q} \} \ \mbox{is stationary}
\end{equation}
for $\ell = 1,2$.
 Then 
\[ p \forcesP \{ \delta \in T_\ell \mid p^\ell_\delta \in G_P\}\ \mbox{is stationary}.\]
Let $G\subset P$ be $V$-generic, with $p \in G$. In $V[G]$ form the 
interpretations $q[G] $ (interpretation of $q$) and $Q[G]$ (interpretation
of $Q$). Then $q[G] \in Q[G]$. 
Define the sets
\[ T'_\ell = \{ \delta \in T_\ell \mid p_\delta^\ell \in G
\},\ \ell = 1,2\]
(which are stationary) and define the sequences
\[ \langle q^\ell_\delta[G] \mid \delta \in T'_\ell \rangle,\
\mbox{ for } \ell = 1,2 .\]
Then in $V[G]$
\[ q[G] \forces_{Q[G]}\ \{ \delta \in T_\ell ' \mid
q^\ell_\delta [G] \in H \}\ \mbox{ is stationary}
\]
where $H$ is the name for the $V[G]$ generic filter over
$Q[G]$. (This follows from (\ref{St}) and since forcing with
$P\ast Q$ is equivalent to the iteration of forcing with $P$
and then with $Q[G]$.)

Since $Q[G]$ satisfies the polarized chain condition
for $T$, there are stationary sets $T_\delta^{''}
\subseteq T'_\delta$ such that:
\begin{quote}
 if $\delta_1\in
T^{''}_1$, $\delta_2\in T_2^{''}$, and $\delta_1 <
\delta_2$, then $q^1_{\delta_1}[G]$ and
$q^2_{\delta_2}[G]$ are compatible in $Q[G]$.
\end{quote}

Back in $V$, let $S_1$ and $S_2$ be $V^P$ names of
$T^{''}_1$ and $ T_2^{''}$ respectively, forced to
have these properties. The following short lemma will be
applied to $S_1$ and to $S_2$.

\begin{lemma}
Suppose that $\langle p_\delta \mid \delta \in T
\rangle$ is a sequence in $P$, $S$ is a name of a
subset of $\omega_1$ and  $p\in P$ a condition such that
\[ p \forces_P\ S \subseteq \{ \alpha \in T \mid
p_\delta \in G \}\ \mbox{and } S\ \mbox{ is stationary in }
\omega_1.\]
 Then there is a stationary subset $T^\ast \subseteq
 T$, and conditions $p^\ast_\delta$ extending both  $p_\delta$ and $ p$ for
 each $\delta \in T^\ast$ such that $p_\delta^\ast \forces
 \delta \in S$.
 \end{lemma}
 {\bf Proof.} Define $T^\ast$ by the condition that
 $\delta \in T^\ast$ iff $\delta \in T$ and there is a
 common extension of $p$ and $p_\delta$ that
 forces $\delta \in S$. We must prove that $T^\ast$ is
 stationary. If $C\subseteq \omega_1$ is any closed
 unbounded set, find $p'\geq p$ and $\delta \in C$
 such that $p'\forces \delta \in S$. Then $\delta \in T$ and 
 $p'\forces\ p_\delta \in G$. Hence $p_\delta \leq p'$ (because $P$ is separative). 
So $\delta \in T^\ast$.\qed

 Apply the lemma to $S_1$ and find a stationary set
 $T^\ast _1\subseteq T_1$ and conditions
 $p_\delta^{\ast 1 }\geq p_\delta^1, p$, for $\delta \in
 T^\ast_1 $ such that
 \[ p_\delta^{\ast 1} \forces \delta \in S_1.\]
 Then (lemma \ref{CCC}) find an extension of $p$, denoted $p^\ast$, such that
 \[ p^{\ast} \forces \{ \delta \in T^\ast_1\mid
 p_\delta^{\ast 1} \in G \} \ \mbox{is stationary}.\]
 Apply the same argument to $S_2$, and find
a stationary set $T^\ast_2\subseteq T_2$ and conditions
$p^{\ast 2}_\delta \geq p_\delta^2, p^\ast$ for $\delta \in T_2^\ast$ 
such that $p^{\ast
2}_\delta \forces \delta \in S_2$. Now $p^{\ast\ast}\geq
p^\ast$ can be found such that
\[ p^{\ast\ast}\forces \{ \delta \in T^\ast_2 \mid p_\delta
^{\ast 2} \in G \} \ \mbox{is stationary}.\]
Since $P$ satisfies the p.c.c., there are stationary sets
$T^{\ast\ast}_1 \subseteq T^\ast_1$ and $T^{\ast\ast}_2
\subseteq T^\ast_2$ such that for every $\delta_1 <
\delta_2$ in $T^{\ast\ast}_1$ and  $T^{\ast\ast}_2$  
(respectively) $p^{\ast 1}_{\delta_1}$ and  $p^{\ast
2}_{\delta_2}$ are compatible in $P$, say by some condition
$p'$ extending both. But then $p' \forces \delta_1\in S_1\
\mbox{and } \delta_2 \in S_2$. It follows that
$(p^1_{\delta_1},q^1_{\delta_1})$ and
$(p^2_{\delta_2},q^2_{\delta_2})$ are compatible in $P\ast Q$
showing that $P\ast Q$ satisfies the p.c.c.
The point is that 
 $$p' \forcesP q^1_{\delta_1}\ \mbox{and}\ q^2_{\delta_2}\ \mbox{are
 compatible in}\ Q$$
 and hence for some $q'\in V^P$, $p'\forcesP\ q' \geq
 q^1_{\delta_1},  q^2_{\delta_2}$. That is, $(p'q') \geq
 (p^1_{\delta_1}, q^1_{\delta_1}),\ (p^2_{\delta_2},
 q^2_{\delta_2})$. \qed

\begin{theorem}
Let $T$ be a stationary subset of $\omega_1$.
An iteration with finite support of p.c.c. for $T$ posets is again
p.c.c. for $T$.
\end{theorem}
{\bf Proof.} The theorem is proved by induction on the length,
$\delta$, of the iteration.
For $\delta$ a successor ordinal, this is essentially Lemma
\ref{Pres}. So we assume that $\delta$ is a limit ordinal, and
$\langle P_\alpha \mid \alpha \leq \delta \rangle$ is a finite
support iteration, where $P_{\alpha+1}= P_\alpha \ast Q_\alpha$
is obtained with $Q_\alpha \in V^{P_\alpha}$ a p.c.c poset for $T$.
Thus conditions in $P_\delta$ are finite functions $p$ defined on a
finite subset $\dom(p)\subset \delta$, and are such that for every
$\alpha \in \dom(p)$, $p\restriction \alpha \forces_{P_\alpha}\
p(\alpha) \in Q_\alpha$. It is well-known that $P_\delta$ satisfies
the c.c.c, and we must prove the polarized property.

Suppose that $\overline{p}^\ell = \langle p_i ^\ell \mid i \in
T_\ell \rangle$ for $\ell = 1,2$ are two sequences of conditions in
$P_\delta$, where $T_\ell \subseteq T$ are stationary, and suppose
also that $p^\ast\in P_\delta$ is such that
 \begin{equation}
 \label{Stat}
 p^\ast \forces_{P_\delta}\ \{ i \in T_\ell \mid p^\ell_i \in G
\}\ \mbox{is stationary for } \ell = 1,2.
\end{equation}
We may assume that $p^\ast$ is compatible with every $p^\ell_i$
(just throw away those conditions that are not).
The case $\cf(\delta)>\omega_1$ is trivial, because the support
of all conditions is bounded by some $\delta' < \delta$ to which induction 
is applied). So there are two
cases to consider.

\noindent
{\bf Case 1:} $\cf(\delta)=\omega$. Let $\langle \delta_n \mid n
\in \omega \rangle$ be an increasing $\omega$-sequence converging
to $\delta$. For every $p_i^\ell$ there is $n\in \omega$ such that
$\dom(p_i^\ell)\subseteq \delta_n$. 
It follows from (\ref{Stat}) that for some specific
$n\in \omega$, for some extension $p^{\ast\ast}\geq p^\ast$
\[ p^{\ast\ast}\forces 
\{ i \in T_\ell \mid p_i^\ell \in P_{\delta_n}\cap G \}\ \mbox{is
stationary for } \ell = 1,2.\]
Now we can apply the inductive assumption to $P_{\delta_n}$.

\noindent
{\bf Case 2:} $\cf(\delta)=\omega_1$. Let $\langle \delta_\alpha
\mid \alpha \in \omega_1\rangle$ be an increasing, continuous, and
cofinal in $\delta$ sequence. Intersecting $T_1$ and $T_2$ with a
suitable closed unbounded set, we may assume that for every $\alpha< \beta$
 $\alpha\in T_1$ and $\beta\in T_2$,  $\dom(p_\alpha^1)\subset \beta$.

We claim that we may without loss of generality assume that,
 for some $\gamma < \delta$, $\dom(p_\alpha^\ell)\cap \delta_\alpha$ for all $\alpha \in T_\ell$. We get this in two steps.

In the first step, find a stationary $T_1'\subseteq T_1$ such that
the sets $\dom(p^1_\alpha) \cap \delta_\alpha$, for $\alpha \in T_1'$, 
are bounded by some $\gamma < \delta$.
 For each $\alpha\in T'_1$ let $p^{1\ast}_\alpha$
be a common extension of $p^1_\alpha$ and $p^\ast$. Then (use lemma \ref{CCC}) 
find an extension $p^{\ast\ast}\geq p^\ast$ such that
\[ p^{\ast\ast}\forces \{ \alpha \in T_1'\mid p_\alpha^1\in G \} \
\mbox{is stationary}.\]
Since $p^{\ast\ast}$ extends $p^\ast$, 
$p^{\ast\ast}\forces \{ i \in T_2\mid p_i^2\in G \}\
\mbox{\em is stationary}$. We can again assume that each $p_i^2$ is compatible
with $p^{\ast\ast}$ and get 
 $T_2'\subseteq T_2$ stationary such that $\dom(p^2_\alpha)\cap \delta_\alpha$
is bounded by some $\gamma' < \delta$ (we rename $\gamma$ to be the maximum of
$\gamma$ and $\gamma'$). 
Rename the stationary sets as $T_1$, $T_2$ and we have our assumption.

Apply induction to $P_\gamma$ and to the conditions $p^\ell_\alpha\restriction \gamma$. This yields two stationary subsets which are as required. \qed
\section{The model}
\begin{theorem}
Assuming the consistency of ZFC, the following property is consistent with
ZFC. There is a stationary co-stationary set $S\subseteq \omega_1$
 such that 
\begin{enumerate}
\item 
For every ladder system $C$ over $S$,
 every gap  contains a $C$-Hausdorff subgap.
\item For every ladder system $H$ over $T = \omega_1\setminus S$ 
there is a gap $g$ with no subgap that is
$H$-Hausdorff.
\end{enumerate}
\end{theorem}

To obtain the required generic extension we assume that $\kappa$ is a cardinal
in $V$ (the ground model) such that $\cf(\kappa) > \omega_1$ and even
 $\kappa^{\aleph_1} = \kappa$. We shall obtain a
generic extension $V[G]$ in which $2^{\aleph_0}=\kappa$ and the two required
properties of the theorem hold. For this we define a finite
support iteration of length $\kappa$, iterating posets $P$ as in section
\ref{GF}, which introduce generic gaps, and posets of the form $Q_{g,C}$, as in
section \ref{Sp}, which are designed to introduce a $C$-Hausdorff subgap to
$g$.

We denote this iteration $\langle P_\alpha \mid \alpha < \kappa \rangle$. So
$P_{\alpha + 1} \simeq P_\alpha \ast R(\alpha)$, where the $\alpha$-th iterand
$R(\alpha)$ is either some $P$ or some $Q_{g,C}$. The rules to
determine $R(\alpha)$ are specified below.
For any limit ordinal $\delta \leq \kappa$, $P_\delta$ is the finite support
iteration of the posets $\langle P_\alpha \mid \alpha < \delta \rangle$.
We define $P_\kappa$ as our final poset, and we shall 
prove that in $V^{P_\kappa}$ the two properties of the theorem hold.

Recall that $P$ satisfies Talayaco's condition and is hence a p.c.c. poset,
and each $Q_{g,C}$ is p.c.c. over $T=\omega_1 \setminus S$ (by lemma
\ref{LCC}).

Since the iterand posets satisfy the p.c.c over $T$, 
each $P_\alpha$ is a p.c.c. poset over $T$
(and in particular a c.c.c poset). It follows that every ladder system and
every gap in $V^{P_\kappa}$ are already in some $V^{P_\alpha}$ for 
$\alpha < \kappa$. It is obvious that if $g\in V^{P_\alpha}$  is 
forced by $p\in P_\kappa$
to be a gap in $V^{P_\kappa}$, then $p\restriction \alpha$ forces it
to be a gap already in $V^{P_\alpha}$. 

To determine the iterands, 
we assume a standard bookkeeping scheme which ensures two things:
\begin{enumerate}
\item
For every ladder system $C$ over $S$ and gap $g$ in $V^{P_\kappa}$, there
exists a stage $\alpha < \kappa$ so that $C,g \in V^{P_\alpha}$, and $R(\alpha)$ is
$Q_{g,C}$.
\item
For some  unbounded  set of ordinals $\alpha \in \kappa$ the 
iterand $R(\alpha)$ is $P$, producing a generic gap $g$, and the subsequent 
iterand $R(\alpha+1)$
is $Q_{g,C}$ for some ladder sequence $C$ over $S$.
\end{enumerate}

The first item ensures that, in $V^{P_\kappa}$, for every ladder system $C$
over $S$,
 every gap  contains a $C$-Hausdorff subgap. (A $C$-Hausdorff subgap in
 $V_{\alpha+1}$ remains $C$-Hausdorff at every later stage and in the
 final model).

 Suppose now that $H$ is a ladder over $T = \omega_1\setminus S$. Then $H$
 appears in some $V^{P_\alpha}$ such that $R(\alpha)$ is the poset $P$, and
 $R(\alpha+1)$ is the poset $Q_{g,C}$ where $g$ is the generic gap
 introduced by $R(\alpha)$, and $C$ is some ladder
 sequence over $S$.
We want to prove that $g$ is
 a gap that has no $H$-Hausdorff subgap in $V^{P_\kappa}$.
 We first prove that $g$ remains a gap in $V^{P_\kappa}$. It is
 clearly a gap in $V^{P_{\alpha+1} }$ by Lemma \ref{Gap}. 
 Since $g$ is $C$-Hausdorff in $V^{P_{\alpha+2} }$, 
 it remains a gap in $V^{P_\kappa}$ (by
 Lemma \ref{StatGap}).

 This generic gap $g$ satisfies the conclusion of lemma
 \ref{Prop} in $V^{P_{\alpha+1} }$:
\begin{quote}
If $J, K \subseteq \omega_1$ are unbounded, then
there is  a club set $D_0\subseteq \omega_1$ such that for every
$\delta \in D_0 $ and
$k\in K\setminus \delta$ there are $m\in \omega$ and a sequence
$j(n)\in \delta \cap J$ increasing and cofinal in $\delta$ such that
 $a_{j(n)}\setminus m \subset b_k$ for all $n\in \omega$.
 \end{quote}

 Since $P_\kappa \simeq P_{\alpha+1} \ast R$, where the remainder $R \simeq
 P_\kappa/P_{\alpha+1}$ is interpreted in $V^{P_{\alpha+1}}$ as a finite support
 iteration of p.c.c. posets over $T$, 
we can view $P_\kappa$ as a two-stage iteration in
 which the second stage is a p.c.c. poset over $T$. 
Thus, for simplicity of expression,
 we can assume that $V^{P_{\alpha+1}}$ is the ground model.
The following lemma then ends the proof.

 \begin{lemma}
 Suppose in the ground model $V$  a ladder
 system $H$ over a stationary set $T\subseteq \omega_1$, and a gap $g$
 that has the property quoted above (the conclusion of lemma \ref{Prop}). 
 Suppose also a poset $R$ that is
 p.c.c. over $T$. Then in $V^R$ the gap $g$ contains no
 $H$-Hausdorff subgap.
 \end{lemma}
{\bf Proof.} Let $g = \{ (a_i \mid i \in \omega_1), (b_j \mid j \in
\omega_1)\}$ and assume (for the sake of a contradiction)
some condition $q'$  in $R$ forces that
$g' = \{ (a_\alpha \mid \alpha \in A), (b_\beta \mid \beta \in B)\}$ 
is a $H$-Hausdorff subgap, where
$A$ and $B$ are names forced by $q$ to be unbounded in $\omega_1$.
Since every club subset of $\omega_1$
in a c.c.c. generic extension contains a club subset in
the ground model, we may assume that the club, $D$, which appears 
in definition \ref{DHG} (of $g'$  being $H$-Hausdorff) is in $V$.

For every $\delta\in T \cap D$ define two conditions in $R$ (extending
the given condition $q$):
\begin{enumerate}
\item
$p_\delta \in R$ is such that 
for some $\alpha(\delta) \in \omega_1\setminus \delta$,
$p_\delta\forcesR \alpha(\delta) \in A$.
(This is possible since $A$ is forced to be unbounded.)
\item
$q_\delta\in R$ extending $p_\delta$ is such that, 
for some $\beta(\delta) \in \omega_1 \setminus \delta$,
 $ q_\delta \forcesR\ \beta(\delta) \in B$. Moreover, as $g'$ is forced to be
$H$-Hausdorff, we can assume that for some 
 $m_\delta \in \omega$,
 $$q_\delta \forcesR\  \mbox{for every}\ n \geq m_\delta,\ \mbox{if}\ i \in A
\cap (\delta \setminus c_\delta(n)),\ \mbox{then}\ X(a_i, b_{\beta(\delta}))> n.$$

\end{enumerate}
By Lemma \ref{CCC} some condition forces that $q_\delta\in G$ (and hence
$p_\delta\in G$) for a stationary set of indices $\delta\in T\cap D$.
Since $R$ is p.c.c. for $T$, there are stationary subsets $T_1, T_2 \subseteq
T$ such that any $p_{\delta_1}$ is compatible with $q_{\delta_2}$ if $\delta_1
\in T_1$, $\delta_2\in T_2$ and $\delta_1 < \delta_2$. 

Consider now the two unbounded sets $J = \{ \alpha(\delta) \mid \delta
\in T_1 \}$, and $K = \{ \beta(\delta) \mid \delta \in T_2\}$. Apply the
conclusion of lemma \ref{Prop} quoted above to $J$ and $K$,
and let $D_0$ be the club set that appears there.
Pick any $\delta \in D \cap T_2 \cap D_0$.
Consider $k = \beta(\delta)$. Then $k \in K\setminus \delta$, and so
there are $m\in \omega$ and a sequence $j(n)\in \delta \cap J$ cofinal in
$\delta$ such that 
\begin{equation}
\label{Contra}
a_{j(n)} \setminus m \subset b_k\ \mbox{ for all } n\in \omega. 
\end{equation}
Yet every $j(n)$ is of the form $\alpha(\delta_n)$ for some
$\delta_n \in T_1 \cap \delta$, and the $\delta_n$'s tend to $\delta$.
So (\ref{Contra}) can be written as
\begin{equation}
\label{Contra1}
X(a_{\alpha(\delta_n)}, b_k) \leq m.
\end{equation}
It  follows from the definition of $T_1$ and $T_2$ 
 that $p_{\delta_n}$ and $q_\delta$ are
compatible in $R$. If $q'$ is a common extension, then $q'$ forces that
for every $n\geq m_\delta$ if $i=\alpha(\delta_n) \geq c_\delta(n)$ then
$X(a_i, b_k) > n$. 
It suffices now to take $n\geq \max \{ m,m_\delta \}$
 to get the contradiction to (\ref{Contra1}).

\end{document}